\documentclass[12pt]{article} 
\usepackage{amstex, amsthm}

\newcommand{\doublespace}{
   \renewcommand{\baselinestretch}{1.2}
   \large\normalsize}

\def \Z{\Bbb Z}
\def \C{\Bbb C}

\def \Q{\Bbb Q}

\def \wt{{\rm wt}}

\def \Res{{\rm Res}}

\def \End{{\rm End}}

\def \<{\langle} 
\def \>{\rangle} 

\def \a{\alpha }

\def \l{\lambda }

\def \b{\beta }
\def \o{\omega }

\newcommand{\1}{\mathbf{1}}
\newcommand{\NO}{\,{\raise0.25em\hbox{$\mathop{\hphantom{\cdot}}%
\limits^{_{\circ}}_{^{\circ}}$}}\,}
\newcommand{\Ker}{\operatorname{Ker}}

\theoremstyle{plain}
 \newtheorem{theorem}{Theorem}[section]
        
        \newtheorem{lemma}[theorem]{Lemma}
        \newtheorem{proposition}[theorem]{Proposition}

 \newtheorem{remark}[theorem]{Remark}

\doublespace

\begin{document}

\newtheorem{thmm}{Theorem}
\newtheorem{thm}{Theorem}[section]
\newtheorem{prop}[thm]{Proposition}
\newtheorem{cor}[thm]{Corollary}
\newtheorem{lem}[thm]{Lemma}
\newtheorem{rem}[thm]{Remark}
\newtheorem{de}[thmm]{Definition}


\begin{center}
{\Large {\bf Classification of irreducible modules for the vertex operator algebra $M(1)^+$ }} \\
\vspace{0.5cm}
Chongying Dong\footnote{Supported by NSF grant 
DMS-9700923 and a research grant from the Committee on Research, UC Santa Cruz.}
and Kiyokazu Nagatomo\footnote{On leave of absence from
Department of Mathematics, Graduate School of Science,
Osaka University, Toyonaka, Osaka 560-0048, Japan.
This work was partly supported by Grant-in-Aid for Scientific Research,
the Ministry of Education, Science and Culture.}
\\
Department of Mathematics, University of California\\
 Santa Cruz, CA 95064\\
\end{center}

\hspace{1.5 cm}
\begin{abstract}
We classify the irreducible modules for 
the fixed point vertex operator subalgebra of the vertex operator algebra
associated to the Heisenberg algebra with the central charge 1
under the $-1$ automorphism.
\end{abstract}

\section{Introduction}
\setcounter{equation}{0}

Let ${\frak h}$ be a finite dimensional complex vector space of dimension 
$d$ with
a nondegenerate symmetric bilinear form and $\hat {\frak h}=
{\frak h}\otimes \C[t,t^{-1}]+\C\,c$ the corresponding affine algebra.
Then the free bosonic 
Fock space $M(1)=S({\frak h}\otimes t^{-1}\C[t^{-1}])$ 
is a vertex operator algebra of central charge $d$ (cf. [FLM]). 
If $d=1$
the automorphism group of $M(1)$ is $\Z_2$ generated by
$\theta$ (see Subsection \ref{subsection:2.3}). Then $M(1)$ has only two proper subalgebras,
namely, $M(1)^+$ and the vertex operator subalgebra 
generated by the Virasoro algebra [DG].
In this paper we determine the Zhu's algebra $A(M(1)^+)$ and
classify the irreducible modules for $M(1)^+.$ The classification
result says that any irreducible $M(1)^+$-module is isomorphic
to either a submodule of a $M(1)$-module  or a submodule of 
$\theta$-twisted $M(1)$-module. 

The vertex operator algebra $M(1)^+$ is closely related to $W$-algebra.
It was shown in [DG] that $M(1)^+$ is generated by the Virasoro element 
$\omega$ and a highest weight vector $J$ of weight 4  
(see Section 2 below). Thus 
$M(1)^+$ can be regarded as the vertex operator algebra associated to
the $W$-algebra $W(2,4)$ with central charge 1 (cf. [BFKNRV]). So this paper
also gives a classification of irreducible modules for the $W$-algebra 
$W(2,4)$  which can be lifted to modules for $M(1)^+.$

The result in this paper is fundamental in the classification of 
irreducible modules for the vertex operator algebras $V_L^+$ [DN].
Let $L$ be a positive definite even lattice of rank 1. The corresponding
vertex operator algebra $V_L$ is a tensor product of $M(1)$ with
the group algebra $\C[L].$ The structure and
representation theory of $V_L$ including the fusion rules
are well understood (see [B], [FLM], [D], [DL], [DLM1]).
Then $\theta$ can be extended to an automorphism of $V_L$ of order 2.
Moreover the fixed point vertex operator
subalgebra $V_L^+$ contains $M(1)^+$ as a subalgebra and $V_L^+$
is a completely reducible $M(1)^+$-module. The result of  the present paper
has been used in [DN] to determine the Zhu's algebra $A(V_L^+)$ and
to classify the irreducible modules for $V_L^+.$

It should be pointed out that conformal field theory associated
to $M(1)^+$ is an orbifold theory (cf. [DVVV]) for the nonrational vertex 
operator algebra $M(1).$ Let $V$ be a rational vertex operator algebra
and $G$ be a finite group of automorphisms of $V.$ The orbifold theory
conjectures that any irreducible module of the fixed point
vertex operator subalgebra $V^G$ is isomorphic to a submodule 
of a $g$-twisted $V$-module for some $g\in G.$ Our result in this
paper suggests that it may be true even when $V$ is not rational. 

One important tool in  the representation theory of vertex
operator algebra is the Zhu's algebra [Z]. In [Z] it was shown
that for any vertex operator algebra $V$ there is an associative
algebra $A(V)$ associated to $V$ such that there is a one to one
correspondence between the irreducible admissible $V$-modules
and irreducible $A(V)$-modules (see Subsection \ref{subsection:2.2}
for more detail).
The main idea in the present paper is to determine the Zhu's algebra
$A(M(1)^+)$ which turns out to be a commutative algebra over $\C$ with
two variables. 

We should mention an important role played by a generalized 
PBW type theorem in this paper. The classical PBW theorem  
gives a basis for the universal enveloping algebra of a Lie algebra
and a nice spanning set of modules. 
For an arbitrary vertex operator algebra $V$ 
the component operators of the fixed generators of $V$ in general
do not form a Lie algebra, so one cannot use the classical PBW theorem  
to get a good spanning set in terms of the component operators of
the generators. As mentioned before $M(1)^+$ is generated by
$\o$ and $J.$ Although the component
operators of $\o$ and $J$ do not form a Lie algebra as the commutators
involves quadratic or higher products, we manage to obtain
a kind of PBW type result
which is good enough to give nice spanning sets of $M(1)^+$ and $A(M(1)^+).$
The same idea and technique have been developed further in [DN]
to yield  nice spanning sets of $V_{L}^+$ and $A(V_L^+).$
A PBW-type generating property 
for general vertex operator algebras has been given in [KL] recently. 

The structure of this paper is as follows. In Section 2 we recall the
definition of admissible twisted modules for a vertex operator algebra,
the notion of Zhu's algebra and related results, and
the construction of vertex operator algebra $M(1)^+.$ In Section 3
we give the commutator relations for the components operators
of $\o$ and $J,$ and we produce
a kind of generalized PBW theorem. This enables us to get spanning
sets of $M(1)^+$ and $A(M(1)^+).$ Section 4 shows how to evaluate
the generators of $A(M(1)^+)$ on the top levels of the known irreducible
modules for $M(1)^+$ to yield the relations which are good enough to
determine the algebra structure of $A(M(1)^+).$ We then 
use $A(M(1)^+)$ to classify the irreducible modules for $M(1)^+.$   

\section{Preliminaries}
\setcounter{equation}{0}

This section is divided into three parts. In the first part we recall
various notions of (twisted) modules for a vertex operator algebra $V$
(cf. [DLM2]). The Zhu's algebra [Z] and related results are explained in 
second part. In the last part we review the vertex operator algebra $M(1)$ 
and its (twisted) modules (cf. [FLM]). 

\subsection{Modules}\label{subsection:2.1}
Let $V$ be a vertex operator algebra (cf. [B], [FLM]) and
$g$ an automorphism of $V$ of finite order
$T.$  Denote the decomposition of $V$ into eigenspaces with respect to
the action of $g$ as $V=\bigoplus_{r\in \Z/T\Z}V^r$
where $V^r=\{v\in V|gv=e^{-2\pi ir/T}v\}$.

An admissible $g$-twisted $V$-module (cf. [DLM2], [Z])  
$$M=\sum_{n=0}^{\infty}M(\frac{n}{T})$$
is an $\frac{1}{T}\Z$-graded vector space with the top level 
$M(0)\ne 0$ equipped 
with a linear map
$$\begin{array}{l}
V\longrightarrow (\End\,M)\{z\}\\
v\longmapsto\displaystyle{ Y_M(v,z)=\sum_{n\in\Q}v_nz^{-n-1}\ \ \ (v_n\in
\End\,M)}
\end{array}$$
which satisfies the following conditions;
for all $0\leq r\leq T-1,$ $u\in V^r$, $v\in V,$ 
$w\in M$,
\begin{eqnarray*}
& &Y_M(u,z)=\sum_{n\in \frac{r}{T}+\Z}u_nz^{-n-1}, \label{1/2}\\ 
& &u_nw=0\ \ \                                  
\mbox{for}\ \ \ n\gg 0,\label{vlw0}\\
& &Y_M({\bf 1},z)=1,\label{vacuum}
\end{eqnarray*}
 \begin{equation}\label{eqn:2.1}
\begin{array}{c}
\displaystyle{z^{-1}_0\delta\left(\frac{z_1-z_2}{z_0}\right)
Y_M(u,z_1)Y_M(v,z_2)-z^{-1}_0\delta\left(\frac{z_2-z_1}{-z_0}\right)
Y_M(v,z_2)Y_M(u,z_1)}\\
\displaystyle{=z_2^{-1}\left(\frac{z_1-z_0}{z_2}\right)^{-r/T}
\delta\left(\frac{z_1-z_0}{z_2}\right)
Y_M(Y(u,z_0)v,z_2)},
\end{array}
\end{equation}
where $\delta(z)=\sum_{n\in\Z}z^n$ (elementary
properties of the $\delta$-function can be found in [FLM]) and
all binomial expressions (here and below) are to be expanded in nonnegative
integral powers of the second variable;
$$u_mM(n)\subset M(\wt(u)-m-1+n)$$
if $u$ is homogeneous. 
If $g=1$, this reduces to the definition of an admissible $V$-module.

A $g$-{\em twisted $V$-module} is
an admissible $g$-twisted $V$-module $M$ which carries a 
$\C$-grading induced by the spectrum of $L(0).$ That is, we have
$$M=\coprod_{\lambda \in{\C}}M_{\lambda} $$
where $M_{\l}=\{w\in M|L(0)w=\l w\}.$ Moreover we require that 
$\dim M_{\l}$ is finite and for fixed $\l,$ $M_{{n\over T}+\l}=0$
for all small enough integers $n.$ Again if $g=1$  we get an ordinary
$V$-module.

\subsection{Zhu's algebra}\label{subsection:2.2}

Let us recall that a vertex algebra $V$ is $\Z$-graded:
\[
V = \coprod_{n\in\Z}V_n,\quad v\in V_n,\quad n = \wt (v).
\]
Each $v\in V_n$ is a homogeneous vector of weight $n$.
In order to define Zhu's algebra $A(V)$ we need two
products $*$ and $\circ$ on $V$. For $u\in V$ homogeneous and $v\in V$ 
\begin{align}\label{eqn:2.2}
u*v&={\rm Res}_z\left(
\frac{(1+z)^{\wt(u)}}{z}Y(u,z)v
\right)
 =\sum_{i=0}^{\infty}{\wt(u)\choose i}u_{i-1}v,\\
\label{eqn:2.3}
u\circ v&={\rm Res}_z\left(
\frac{(1+z)^{\wt(u)}}{z^2}Y(u,z)v\right) = \sum_{i=0}^\infty\binom{\wt(u)}{i}
u_{i-2}v
\end{align}
and extend  both (\ref{eqn:2.2}) and (\ref{eqn:2.3}) 
to linear products on $V.$ Define $O(V)$ to
be the linear span of all $u\circ v$ for $u,v\in V.$ 
Set $A(V)=V/O(V).$ 
For $u\in V$ we denote $o(u)$ the weight zero 
component operator of $u$ on any admissible module. Then $o(u)=u_{\wt(u)-1}$
if $u$ is homogeneous. The following theorem is essentially due to Zhu [Z].
 
\begin{theorem} 
\label{theorem:2.1} 
{\rm (i)} The product $*$ induces 
an associative algebra structure on $A(V)$ 
with the identity $\1+O(V).$ Moreover $\o+O(V)$ is a central
element of $A(V).$ 

{\rm (ii)} The map $u\mapsto o(u)$ gives a representation of $A(V)$ on $M(0)$
for any admissible $V$-module $M.$ Moreover, if any admissible
$V$-module is completely reducible, then $A(V)$ is a finite dimensional
semisimple algebra. 

{\rm (iii)} The map $M\to M(0)$ gives a bijection between 
the set of equivalence classes of simple admissible $V$-modules and the 
set of equivalence classes of simple $A(V)$-modules.
\end{theorem}

For convenience we write $[u]=u+O(V)\in A(V).$ 
We define $u\sim v$ for $u,v\in V$ if $[u]=[v].$ This induces a relation on 
$\End\, V$ such that for $f,g\in \End V$, $f\sim g$ if and only if
$fu\sim gu$ for all $u\in V$.

The following proposition is useful later (cf. [W], [Z]).
\begin{proposition}\label{proposition:2.2}
{\rm (i)} Assume that $u\in V$ homogeneous, $v\in V$ and $ n \geq 0$.
Then 
$$\Res_z\left(
\frac{(1+z)^{\wt(u)}}{z^{2+n}}Y(u,z)v
\right)
= \sum_{i=1}^\infty \binom{\wt(u)}{i}u_{i-n-2}v
\in O(V).$$
{\rm (ii)} If $u$ and $v$ are homogeneous elements of $V$, then
$$u*v\sim
\Res_z\left(
\frac{(1+z)^{\wt(v)-1}}{z}Y(v,z)u
\right).$$
{\rm (iii)} For any $n\geq 1$, 
\begin{equation}\label{eqn:2.6}
L(-n)\sim (-1)^n
\left\{
(n-1)(L(-2)+L(-1))+L(0)
\right\}
\end{equation}
where $L(n)$ are the Virasoro operators given by $Y(\o,z)=\sum_{n\in \Z}L(n)z^{-n-2}.$ 
\end{proposition}

\subsection{Vertex operator algebras $M(1)$ and $M(1)^+$}\label{subsection:2.3}

Finally we discuss the construction of vertex operator algebra $M(1)$
and its (twisted) modules (cf. [FLM]). We also define the 
vertex operator subalgebra $M(1)^+.$  

Let ${\frak h}$ be a finite dimensional vector space
with a nondegenerate symmetric bilinear form 
$\<\cdot,\cdot\>$ and $\hat {\frak h}={\frak h}\otimes 
\C[t,t^{-1}]\oplus \C\,c$ the corresponding affine Lie algebra.
Let $\l\in \frak{h}$ and consider the
induced $\hat {\frak h}$-module
$$M(1,\l)= 
U(\hat{\frak h})\otimes_{U({\frak h}\otimes{{\C}}[t]
\oplus{{\C}}c)}{{\C}}\simeq S(\frak h\otimes t^{-1}\C[t^{-1}])\ \ \ (\mbox{linearly})
$$
where ${\frak h}\otimes t{{\C}}[t]$ acts trivially on $\C,$
${\frak h}$ acts as $\<\alpha,\l\>$ for $\alpha\in{\frak h}$ and $c$
acts as 1.
For $\alpha\in{\frak h}$ and $n\in {\Z}$,  we write $\alpha(n)$
for the operator $\alpha\otimes t^n$ and put
$$\alpha(z)=\sum_{n\in{{\Z}}}\alpha(n)z^{-n-1}.$$
Among $M(1,\lambda),\lambda \in \mathfrak{h}$, 
$M(1) = M(1,0)$ has the special interesting 
because it has a natural vertex operator algebra structure
as explained below. For $\alpha_1, ..., \alpha_k \in
{\frak h},\  n_1, ..., n_k \in{\Z}\, (n_i>0)$ and  $v=\alpha_1(-n_1)\cdot
\cdot\cdot\alpha_k(-n_k)\in M(1)$, we define a vertex operator 
corresponding to $v$ by
$$Y(v,z)=\mbox{$\circ\atop\circ$}
\partial_z^{(n_1-1)}\alpha_1(z)\partial_z^{(n_2-1)}\alpha_2(z)\cdots
\partial_z^{(n_k-1)}\alpha_k(z)\mbox{$\circ\atop\circ$},
$$
where
\[
\partial_z^{(n)} = \frac{1}{n!}
\left(
\frac{d}{dz}
\right)^n
\]
and a normal ordering procedure indicated by open colons signifies
that the expression above is to be reordered if necessary so that
all the operators
$\alpha(n)$ $(\alpha\in{\frak h},\  n<0)$  are to be placed to
the left of all the operators $\alpha(n)$ $(n\ge 0)$ before the expression 
is evaluated. We extend $Y$ to all $v\in V$ by linearity. 
Let $\{\b_1,...\b_d\}$ be
an orthonomal basis of ${\frak h}.$ Set ${\bf 1}=1$ and 
$\omega=\frac{1}{2}\sum_{i=1}^d\b_i(-1)^2$ The following theorem is well known
(cf. [FLM]):
\begin{theorem}\label{theorem:2.4}
The space $M(1)=(M(1),Y,{\bf 1},\omega)$ is a simple vertex 
operator algebra and $M(1,\l)$ for $\l\in {\frak h}$ gives a complete
list of inequivalent irreducible modules for $M(1).$
\end{theorem}

We define an automorphism $\theta$ of $M(1)$ by 
$$\theta(\a_1(n_1)\cdots
\a_k(n_k))=(-1)^k\a_1(n_1)\cdots\a_k(n_k).$$ Then 
$\theta$-invariants $M(1)^+$ of $M(1)$ form a simple vertex operator 
subalgebra and the $-1$-eigenspace $M(1)^-$ is an irreducible $M(1)^+$-module
(see Theorem~2 of [DM2]). Clearly $M(1)=M(1)^+\oplus M(1)^-$. 

Following [DM1] we define $\theta\circ M(1,\l)=(\theta\circ M(1,\l),
Y_{\theta})$ where $Y_{\theta}(v,z) =Y(\theta v,z).$ Then $\theta\circ
M(1,\l)$ is also an irreducible $M(1)$-module isomorphic to $M(1,-\l).$
The following proposition is an direct consequence of Theorem 6.1 of [DM2]:

\begin{proposition}\label{proposition:2.5} 
If $\l\ne 0$ then $M(1,\l)$ and $M(1,-\l)$ are
isomorphic and irreducible $M(1)^+$-modules.
\end{proposition}

Next we turn our attention the $\theta$-twisted $M(1)$-modules (cf. [FLM]).
The twisted affine algebra is defined to be $\hat {\frak h}[-1]=
\sum_{n\in\Z}{\frak h}\otimes t^{1/2+n}\oplus \C\, c$ and its canonical 
irreducible module is 
$$M(1)(\theta)=U(\hat {\frak h}[-1])\otimes_{U({\frak h}\otimes t^{1/2}\C[t]\oplus \C 
c)}\C\simeq S({\frak h}\otimes t^{-1/2}\C[t^{-1/2}])$$
where ${\frak h}\otimes t^{1/2}\C[t]$ acts trivially on $\C$ and
$c$ acts as 1. As before 
there is an action of $\theta$ on $M(1)(\theta)$ by $\theta (\a_1(n_1)\cdots \a_k(n_k))
=(-1)^k\a_1(n_1)\cdots \a_k(n_k)$ where $\a_i\in {\frak h},$
$n_i\in \frac{1}{2}+\Z$ and $\alpha(n)=\alpha\otimes t^n$. We denote
the $\pm 1$-eigenspace of $M(1)(\theta)$ under $\theta$ 
by $M(1)(\theta)^{\pm}.$

Let $v=\a_1(-n_1)\cdots \a_k(-n_k)\in M(1)$. We define 
$$W_{\theta}(v,z)=\mbox{$\circ\atop\circ$}
\partial_z^{(n_1-1)}\alpha_1(z)\partial_z^{(n_2-1)}\alpha_2(z)
\cdots\partial_z^{(n_k-1)}
\alpha_k(z)\mbox{$\circ\atop\circ$},
$$
where the right side is an operator on $M(1)(\theta)$, namely,
\[
\alpha(z)=\sum_{n\in \frac{1}{2}+\Z}\alpha(n)z^{-n-1}
\]
and the normal ordering notation
is obvious. Further we extend this to all $v\in V_L$ by linearity.  
Define constants $c_{mn}\in \Q$
for $m, n\ge0$ by the formula
$$\sum_{m,n\ge 0}c_{mn}x^my^n=-{\rm log}\left(\frac
{(1+x)^{1/2}+(1+y)^{1/2}}{2}\right).$$    
Set 
$$\Delta_z=\sum_{m,n\ge 0}\displaystyle{
\sum^d_{i=1}}
c_{mn}\beta_i(m)\beta_i(n)z^{-m-n}.$$
Now we define $twisted$ $vertex$ $operators$
$Y_{\theta}(v,z)$ for $v\in M(1)$
as follows:
$$Y_{\theta}(v,z)=W_{\theta}(e^{\Delta_z}v,z).
$$
 
Then we have 
\begin{theorem}\label{theorem:2.6} {\rm (i)} $(M(1)(\theta),Y_{\theta})$ is 
the irreducible
$\theta$-twisted $M(1)$-module.

{\rm (ii)} $M(1)(\theta)^{\pm}$ are irreducible $M(1)^+$-modules.
\end{theorem}

Part (i) was a result of Chapter 9 of [FLM] and part (ii) follows
Theorem 5.5 of [DLi].

In this paper, we mainly consider the case that $\mathfrak{h}$ is 
one dimensional. From now on we always assume that 
$\mathfrak{h} = \C\,h$ with the normalized inner product
$\langle h,h\rangle = 1$.

\begin{remark}{\rm It is easy to see in this case that the automorphism
group of $M(1)$ is generated by $\theta.$ It was pointed out in [DG] that
$M(1)^+$ is the only proper vertex operator subalgebra of $M(1)$ which 
differs from the vertex operator subalgebra generated by $\o.$  }
\end{remark}

For the later purpose we need to know the 
first few coefficients of $z$ in $\Delta_z.$ Note that 
$$
\begin{aligned}
-{\rm log}&\left(\frac
{(1+x)^{1/2}+(1+y)^{1/2}}{2}\right)\\
=& -\frac{1}{4}x -\frac{1}{4}y 
+\frac{3}{32}x^2 + \frac{1}{16}xy + \frac{3}{32}y^2\\
&-\frac{5}{96}x^3- \frac{1}{32}x^2y -\frac{1}{32}xy^2 -\frac{5}{96}y^3\\
&+\frac{35}{1024}x^4 + \frac{5}{256}x^3y+\frac{9}{512}x^2y^2
+ \frac{5}{256}xy^3 +\frac{35}{1024}y^4 +\cdots. 
\end{aligned}
$$
Thus
\begin{equation}\label{eqn:2.13}
\begin{aligned}
\Delta_z = &-\frac{1}{2}h(0)h(1)z^{-1} +
\left(\frac{3}{16}h(0)h(2)+\frac{1}{16}h(1)^2\right)z^{-2}\\
&+\left(-\frac{5}{48}h(0)h(3)-\frac{1}{16}h(1)h(2)\right)z^{-3}\\
&+\left(\frac{35}{512}h(0)h(4)+\frac{5}{128}h(1)h(3)+\frac{9}{512}h(2)^2
\right)z^{-4}
+\cdots.
\end{aligned}
\end{equation}

\section{A spanning set of $A(M(1)^+)$}
\setcounter{equation}{0}

In this section we use a result in [DG] to yield a spanning set 
of $M(1)^+$  and then use it to produce a spanning
set of $A(M(1)^+).$ We also list known irreducible modules for $M(1)^+$
and the actions of $L(0)$ and $o(J)$ on the top levels of these
modules where $J$ is a singular vector of $M(1)^+$ of weight 4
defined in Subsection \ref{subsection:3.1}.  

\subsection{Some commutator relations}
\label{subsection:3.1}

Recall that $Y(\o,z)=\sum_{n\in \Z}L(n)z^{-n-2}$ 
where the component operators $L(n)$ together with $1$ spanned 
a Virasoro algebra of central charge 1 on $M(1).$
It is well known that $M(1)$ is a unitary representation for the Virasoro
algebra and $M(1)^+,$ as the submodule for the Virasoro algebra, 
is
a direct sum of irreducible modules
\begin{equation}\label{eqn:3.1}
M(1)^+ = \bigoplus_{m\in\Z_{\geq 0}}L(1,4m^2)
\end{equation}
where $L(1,4m^2)$ is an irreducible highest weight Virasoro module with 
highest weight $4m^2$ and central charge $1$ 
(See [DG] Theorem 2.7 (1)).

Let 
\begin{equation}
J = h(-1)^4\1 -2h(-3)h(-1)\1 + \frac{3}{2}h(-2)^2\1\label{eqn:3.2}
\end{equation}
which is a singular vector of weight $4$ for the Virasoro algebra.
Then the field
$$J(z) =\NO h(z)^4\NO-\NO\partial^2_z h(z)h(z)\NO
+\frac{3}{2}\NO(\partial_z h(z))^2\NO
$$
is a primary field. We have commutation relations
$$[L(m),J(z)] = z^{m}(z\partial_z + 4(m+1))J(z)\quad(m\in\Z)
$$
which follows from the Jacobi identity (\ref{eqn:2.1}) and 
which is equivalent to 
\begin{equation}\label{eqn:3.5}
[L(m),J_n] = (3(m+1)-n)J_{n+m}\quad (m,n\in\Z)
\end{equation}
where $J(z) = \sum_{n\in\Z}J_nz^{-n-1}$.

Next we compute the commutator relation $[J_m,J_n]$ for $m,n\in\Z.$
Again by the Jacobi identity (\ref{eqn:2.1}) we know 

\[
[J_m,J_n] = \sum_{i=0}^\infty\binom{m}{i}(J_iJ)_{m+n-i}.
\]
Since the weight of $J$ is 4, we see that $\wt(J_iJ) = 7-i\leq 7.$ 
Then it follows from the decomposition (\ref{eqn:3.1}) that 
for any $i\in\Z_{\geq 0}$, we have
$J_iJ\in L(1,0)\bigoplus L(1,4)$
and then all these are expressed as linear combinations of
\[
L(-m_1)\cdots L(-m_s)\1,\quad
L(-n_1)\cdots L(-n_t)J
\]
where $m_1\geq m_2\geq\dots\geq m_2\geq 2$, $n_1\geq n_2\geq \dots\geq n_t\geq 1$
and $s,t\leq 3$. Note that for any vertex operator algebra $V,$ $u,v\in V$
and $m,n\in\Z,$ $(u_mv)_n$ is a linear combination of operators 
$u_sv_t$ and $v_tu_s$ for $s,t\in \Z.$  Using (\ref{eqn:3.5}) we obtain
the following lemma.
\begin{lemma}\label{lemma:3.1}
For any $m,n\in\Z$, commutators $[J_m,J_n]$ are expressed as linear combinations of
\[
L(p_1)\cdots L(p_s),\quad
L(q_1)\cdots L(q_t)J_r
\]
where $p_1,\dots,p_s,q_1,\dots,q_t,r\in\Z$ and $s,t\leq 3$.
\end{lemma}

\subsection{A spanning set for $M(1)^+$ }
\label{subsection:3.2}

We first note the following theorem.

\begin{theorem}{\rm ([DG], Theorem 2.7 (2))}\label{theorem:3.2}
As a vertex operator algebra, $M(1)^+$ is generated by the Virasoro element 
$\omega$ and any singular vector of weight greater than 0. In
particular, $M(1)^+$ is generated by $\o$ and $J.$
\end{theorem}

{}From this theorem we see that $M(1)^+$ is spanned by
$$\{u^1_{m_1}\cdots u^k_{m_k}\1|u^i=\o, J, m_i\in \Z\}$$
which are not necessarily linearly independent. 
We say that an expression $u^1_{m_1}\cdots u^k_{m_k}\1$ has length $t$ 
with respect to $J$,  which we write $\ell_J(u^1_{m_1}\cdots u^k_{m_k}\1) = t,$
if $\{i|u^i=J\}$ has cardinality $t.$ Note that $\o_i=L(i-1).$ 
An induction on $\ell_J(u^1_{m_1}\cdots u^k_{m_k}\1)$ using (\ref{eqn:3.5}) and 
Lemma \ref{lemma:3.1} shows that $u^1_{m_1}\cdots u^k_{m_k}\1$ is a linear combination of vectors 
of type
\begin{align*}
&\left\{
L(m_1)L(m_2)\cdots L(m_s)J_{n_1}J_{n_2}\cdots J_{n_t}\1\,|\,
m_a,n_b\in\Z
\right\}.
\end{align*}
Thus $M(1)^+$ is spanned by those vectors.

Using the commutator relations (\ref{eqn:3.5}) and 
fact that $L(m)\1= 0, m\geq -1$, we get the following lemma.

\begin{lemma}\label{lemma:3.3}
Let \,$W$ be a subspace of $M(1)^+$ spanned by $J_{n_1}\cdots J_{n_t}\1$
with $n_i\in \Z.$ Then $W$ is invariant under the action of $L(m), m\geq -1$.
\end{lemma}

\begin{proposition}\label{proposition:3.4}
The vertex operator algebra $M(1)^+$ is spanned by
\[
L(-m_1)\cdots L(-m_s)J_{-n_1}\cdots J_{-n_t}\1
\]
where $m_1\geq m_2\geq \dots\geq m_s\geq 2,\, n_1\geq n_2\geq\dots\geq
 n_t\geq 1$.
\end{proposition}
\begin{proof} We have already known that $M(1)^+$ is spanned by
\[
L(-m_1)\cdots L(-m_s)J_{-n_1}\cdots J_{-n_t}\1
\]
where $m_a,n_b\in\Z$. Using the PBW theorem for the Virasoro algebra
we can assume that $m_1\geq \cdots \geq m_s.$ By 
Lemma \ref{lemma:3.3} we can further assume that
$m_1\geq m_2\geq \dots\geq m_s\geq 2$.  
We proceed by induction on the length with respect
to $J$ that $v=L(-m_1)\cdots L(-m_s)J_{-n_1}\cdots J_{-n_t}\1$
can be spanned by the indicated vectors in the proposition.

If the length is $0$, it is clear.
Suppose that it is true for all monomials $v$ such that 
$\ell_J(v)<t$. Since $J_k\1=0$ for $k\geq 0$ we can assume $n_t\geq 1.$
If $n_1\geq \cdots \geq n_t$ we are done. Otherwise 
there exists $n_a$ such that $n_{a+1}\geq \cdots \geq n_t$ but
$n_a< n_{a+1}.$ There are two cases $n_a\leq 0$ and $n_a >0$ 
which are dealt with separately.
If $n_a\leq 0$ then $J_{-n_a}\1 = 0$ and  
\begin{multline*}
L(-m_1)\cdots L(-m_s)J_{-n_1}\cdots J_{-n_t}\1\\
= \sum_{j=a+1}^t L(-m_1)\cdots L(-m_s)J_{-n_1}\cdots
\overset{\,\,\,\vee}{J}_{-n_a}
\cdots [J_{-n_a},J_{-n_j}]\cdots J_{-n_t}\1
\end{multline*}
where $\overset{\,\,\,\vee}{J}_{-n_a}$ means that we omit the term 
$J_{-n_a}$.
However by Lemma \ref{lemma:3.1}, $[J_{-n_a}, J_{-n_j}]$ are 
linear combinations of operators of type
$$L(p_1)\cdots L(p_{s'}),\quad L(q_1)\cdots L(q_{t'})J_r.$$
 By substituting these into the above and using commutation relation
(\ref{eqn:3.5}) again, the right hand side is a linear combination 
of monomials whose lengths with respect to $J$ are 
less than or equal to $t-1$. Thus by induction hypothesis,
this is expressed as linear combinations of expected monomials. 

If $n_a>0$ then either $n_a<n_t$ or there exists $b$ with $t>b>a$ so that
$n_b>n_a \geq n_{b+1}.$ Then we have either  
\begin{align*}
L(-m_1)\cdots& L(-m_s)J_{-n_1}\cdots J_{-n_t}\1\\
&= \sum_{j=a+1}^t L(-m_1)\cdots L(-m_s)J_{-n_1}\cdots
\overset{\,\,\,\vee}{J}_{-n_a}
\cdots [J_{-n_a},J_{-n_j}]\cdots J_{-n_t}\1\\
&\quad\quad\quad+ L(-m_1)\cdots L(-m_s)J_{-n_1}\cdots
\overset{\,\,\,\vee}{J}_{-n_a}
\cdots J_{-n_t}J_{-n_a}\1
\end{align*}
or
\begin{align*}
L(-m_1)\cdots& L(-m_s)J_{-n_1}\cdots J_{-n_t}\1\\
&= \sum_{j=a+1}^b L(-m_1)\cdots L(-m_s)J_{-n_1}\cdots
\overset{\,\,\,\vee}{J}_{-n_a}
\cdots [J_{-n_a},J_{-n_j}]\cdots J_{-n_t}\1\\
&\quad\quad+ L(-m_1)\cdots L(-m_s)J_{-n_1}\cdots
\overset{\,\,\,\vee}{J}_{-n_a}
\cdots J_{-n_b}J_{-n_a}J_{-n_{b+1}}\cdots J_{-n_t}\1
\end{align*}
{}From the discussion of case $n_a\leq 0$ it is enough to show
either 
$$ L(-m_1)\cdots L(-m_s)J_{-n_1}\cdots\overset{\,\,\,\vee}{J}_{-n_a}
\cdots J_{-n_t}J_{-n_a}\1$$
or 
$$ L(-m_1)\cdots L(-m_s)J_{-n_1}\cdots\overset{\,\,\,\vee}{J}_{-n_a}
\cdots J_{-n_b}J_{-n_a}J_{-n_{b+1}}\cdots J_{-n_t}\1$$
can be expressed as linear combinations of desired vectors. But this
follows from an induction on $a.$ 
\end{proof}

\subsection{A spanning set for $A(M(1)^+)$}
\label{subsection:3.3}

For short we set
\[
v^{*s} =\stackrel{s}{\overbrace{v*\cdots *v}},
\]
for $v\in M(1)^+.$ Recalling $[v]=v+O(M(1)^+)$ for $v\in M(1)^+$.
We will also use a similar notation $[v]^{\*s}.$ 
Then it is easy to see that $[v^{*t}]=[v]^{*t}.$

\begin{theorem}\label{theorem:3.5}
The Zhu's algebra $A(M(1)^+)$ is spanned by 
$\mathcal{S} = \{[\omega]^{*s}*[J]^{*t}\,|\, s,t\geq 0\}.$
\end{theorem}
\begin{proof}
By Proposition \ref{proposition:3.4}, 
it is enough to show that for any monomial
\[
v = L(-m_1)\cdots L(-m_s)J_{-n_1}\cdots J_{-n_t}\1
\]
where 
$m_1\geq m_2\geq \dots \geq m_1\geq 2,\,n_1\geq n_2\geq \dots\geq n_t\geq 1,$
$[v]$ is a linear combination of $\mathcal{S}.$
We  prove by induction on $\ell_J(v)$ that $[v]$ is spanned 
by vectors $[\omega]^{*p}*[J]^{*q}$ in  $\mathcal{S}$ such that
$q\leq t$ and  weights of its homogeneous components are less than or
equal to the weight of $v.$   

In case that $\ell_J(v)=0$, then $v = L(-m_1)\cdots L(-m_s)\1$, 
which is spanned by
$\{[\omega]^{*s}|s\geq 0\}$ (cf. [FZ]).
Now let $t>0$ and assume that the statement is true 
for all $v$ with $\ell_J(v)<t$. We will prove
by induction on weight of $v$ that $[v]$ is a linear combination
of $\mathcal{S}.$ 
Clearly, the smallest weight
is $t\,\wt (J)$ and corresponding $v$ has the form 
\[
v = \overset{t}{\overbrace{J_{-1}\cdots J_{-1}}}\1.
\]
Then by (\ref{eqn:2.2}),
\[
J^{*t}-v =
\sum
\begin{Sb}
n_i\in\{-1,0,1,2,3\}\\
(n_i)\neq (-1,\dots,-1)
\end{Sb}
a_{n_1n_2\dots n_t}J_{n_1}J_{n_2}\cdots J_{n_t}\1.
\]
Since each term appeared in the right hand side involves 
$J_{n_i}$ for some nonnegative $n_i$, we can write the right hand side as 
a linear combination of spanning vectors in Proposition \ref{proposition:3.4}
whose lengths are strictly less than $t$. Thus by induction hypothesis,
the image of right hand side in $A(M(1)^+)$ is spanned by $\mathcal{S}$ and
so is $[v].$ 

Now consider general $v = L(-m_1)\cdots L(-m_s)J_{-n_1}\cdots J_{-n_t}\1.$
Without loss of generality, we can assume that $m_1=m_2=\dots=m_s=2$,
namely,
\[
v = \overset{s}{\overbrace{L(-2)\cdots L(-2)}}J_{-n_1}\cdots J_{-n_t}\1
\]
since suppose there exits $m_i$ such that $m_i\geq 3$, then $m_1\geq 3$
and by (\ref{eqn:2.6}), 
\[
v\sim 
(-1)^{m_1}\left\{
(m_1-1)(L(-2)+L(-1))+L(0)
\right\}
L(-m_2)\cdots L(-m_s)J_{-n_1}\cdots J_{-n_t}\1
\]
which is a sum of three homogeneous vectors of weight 
strictly less than $\wt(v)$.
Then we see
\[
v = \omega^{*s}*(J_{-n_1}\1)*
(J_{-n_2}\cdots J_{-n_t}\1) + v'
\]
where $\wt(v')<\wt(v)$. Then again by using induction hypothesis about
weight, it is enough to show that the image of 
\[
v =  \omega^{*s}*(J_{-n_1}\1)*
(J_{-n_2}\cdots J_{-n_t}\1)
\]
in $A(M(1)^+)$ is spanned by $\mathcal{S}$. Since $\omega$ is a central element in 
$A(M(1)^+)$, we have
\begin{align*}
v &= (J_{-n_1}\1)*\omega^{*s}*
(J_{-n_2}\cdots J_{-n_t}\1)\\
&= J_{-n_1}(\omega^{*s}*
(J_{-n_2}\cdots J_{-n_t}\1))+ v'
\end{align*}
where $\wt(v')<\wt(v)$. 
If $n_1>1$, we can use the fact that $J_{-n_1}u$ is congruent
to a sum of vectors whose lengths are less than or equal to $t$ and
whose weights are smaller than $\wt(v)$ (cf. (\ref{eqn:2.3}))
to show that $[v]$ is spanned by $\mathcal{S}.$
If $n_1= 1$, then 
$n_2=\dots=n_t=1$ and
\[
v = \omega^{*s}*J^{*t} 
+\text{lower weight terms}.
\]
Again it is done by the induction assumption. 
\end{proof}

\begin{remark}\label{remark:3.6}
{\rm 
{}From the proof of Theorem \ref{theorem:3.5}, we see 
$v$ is spanned by
$\omega^{*s}*J^{*t}$ with $2s+4t\leq \wt(v)$.
}
\end{remark}

\subsection{List of irreducible modules}
\label{subsection:3.4}

As mentioned in Subsection \ref{subsection:2.3}, $M(1)^+$
has the following irreducible modules,
$$
 M(1)^+,\, M(1)^-,\, M(1,\lambda)\,(0\neq\lambda\in\C),\,
 M(1)(\theta)^+,\, M(1)(\theta)^-.$$
Recall that $M(1,\l)$ and $M(1)(\theta)$ are symmetric algebras
on ${\frak h}\otimes t^{-1}\C[t^{-1}]$ and  
${\frak h}\otimes t^{-1/2}\C[t^{-1}]$ as vectors spaces.  

The following table gives the action of $\omega$ and $J$ on the top levels 
of these modules.

\bigskip
\begin{center}
\begin{tabular}{|c|c|c|c|c|c|}
\hline
&$M(1)^+$&$M(1)^-$&$M(1,\lambda),\,\lambda\in\C^\times$&$M(1)(\theta)^+$&$M(1)(\theta)^-$\\
\hline
$M(0)$&$\C\1$&$\C h(-1)\1$&$\C$&$\C$
&$\C h(-1/2)$\\
\hline
$\omega$&0&1&$\lambda^2/2$&$1/16$&$9/16$\\
\hline
$J$&0&-6&$\lambda^4-\lambda^2/2$&$3/128$&$-45/128$\\
\hline
\end{tabular}
\end{center}

\bigskip
Here we give some explanations on how to get the table.
The actions of $\omega$ and $J$ 
on these spaces except $M(1)(\theta)^\pm$ are easily 
verified. From the definition of $Y_\theta(u,z)$, we see
\[
Y_\theta(\omega,z) = \frac{1}{2}\NO h(z)^2\NO+\frac{1}{16}z^{-2}.
\]
Recall the expression of $J$ from (\ref{eqn:3.2}).
Then by using (\ref{eqn:2.13}), we get
\[
e^{\Delta_z}J =
J + \frac{3}{4}h(-1)^2\1 z^{-2}
+ \frac{3}{128}z^{-4}
\]
and thus
\[
Y_\theta(J,z) = J(z) + \frac{3}{4}\NO h(z)^2\NO z^{-2}+ 
\frac{3}{128}z^{-4}
\]
where $h(z)= \sum_{n\in\frac{1}{2}+\Z}h(n)z^{-n-1}$.
The actions of $\o$ and $J$ on the top levels of $M(1)(\theta)^+$
and $M(1)(\theta)^-$ are immediately derived.

\section{Classification of irreducible modules for $M(1)^+$}
\setcounter{equation}{0}

In this section we explicitly determine the algebra structure of $A(M(1)^+)$
and use this result together to prove that the list of
the irreducible modules in Subsection \ref{subsection:3.4} is complete. 

\subsection{The structure of $A(M(1)^+)$}
\label{subsection:4.1}

It has been proved in Subsection \ref{subsection:3.3} that Zhu's algebra
$A(M(1)^+)$ as an associative algebra is generated by $[\o]$ and $[J].$ 
Since  $[\o]$ is a central element, $A(M(1)^+)$ is a commutative 
associative algebra and  must be isomorphic to 
a quotient of the polynomial algebra $\C[x,y]$ with variables $x$ and $y$
modulo an ideal $I.$ We still need to determine the ideal explicitly.
For this purpose we will find relations between $[\omega]$ and $[J]$ in 
$A(M(1)^+)$. 

For convenience we simply write $[u]$ by $u$ for $u\in M(1)^+$ and
$u*v$  by $uv.$  

\begin{proposition}
\label{proposition:4.1}
In $A(M(1)^+)$ 
\[
J^2 = p(\omega)+q(\omega)J
\]
where
\[
p(x) = \frac{1816}{35}x^4-\frac{212}{5}x^3 +\frac{89}{10}x^2-\frac{27}{70}x,
\quad
q(x)=-\frac{314}{35}x^2+\frac{89}{14}x-\frac{27}{70}.
\]
Or equivalently,
$$(J+\o-4\o^2)(70J+908\o^2-515\o+27)=0.$$
\end{proposition}
\begin{proof} Recall that as a module for the Virasoro algebra $M(1)^+$
has the decomposition
$M(1)^+ = \bigoplus_{m\geq 0}L(1,4m^2)$.
Since $J$ is the singular vector with weight $4$, we see
\[
J^2=\sum_{i\geq 0}{4\choose i}J_{i-1}J\in L(1,0)\bigoplus L(1,4).
\]
Therefore from Remark \ref{remark:3.6}, we get
\begin{equation}\label{eqn:4.1}
J^2 = p(\omega)+q(\omega)J
\end{equation}
where $p,q$ are polynomials of degrees less than or equal to
 $4$ and $2$ respectively.
Let 
 \[
p(x) = \alpha x^4 +\beta x^3 +\gamma x^2 +\delta x +\epsilon
\quad \text{and}\quad
q(x) = ax^2 + bx + c.
\]
In order to determine the coefficients of  $p(x)$ and $q(x)$ 
we evaluate both sides of equation (\ref{eqn:4.1}) on modules 
listed in Subsection \ref{subsection:3.4}.

Since $\omega=J=0$ on the top level of $M(1)^+$ we have $\epsilon = 0$.
On the top level of $M(1)^-,$   $\omega= 1$ and $J =-6$  give
$\alpha+\beta+\gamma+\delta-6(a+b+c)=36$.
Further on the top levels of $M(1,\lambda)$ for $\l\in \C^\times$ 
we know $\omega= \lambda^2/2$ and
$J = \lambda^4-\lambda^2/2$. Comparing the coefficients of $\l^i$'s tells
us 
\[
\alpha+4a= 16,
\quad\beta+4b-a= -8,
\quad\gamma+4c-b=1,
\quad\delta-c=0.
\]
Finally, we get two more equations by substituting 
$\omega= 1/16, J = 3/128$ on $M(1)(\theta)^+$ and
$\omega= 9/16,J = -45/128$ on $M(1)(\theta)^-$. 
Solving this linear system gives the desired result.
\end{proof}

\begin{proposition}
\label{proposition:4.2}
In $A(M(1)^+)$
\[
(\omega-1)(\omega-\frac{1}{16})(\omega-\frac{9}{16})
(J+\omega-4\omega^4) = 0.
\]
\end{proposition}

As a vertex operator algebra, $M(1)^+$ has the weight space decomposition 
$M(1)^+ = \bigoplus_{m\geq 0}M(1)^+_m$. The following is the list of $\dim 
M(1)^+_m$ for $m$ up to 10. 

\bigskip
\begin{center}
\begin{tabular}{|c|c|c|c|c|c|c|c|c|c|c|c|}
\hline
m&0&1&2&3&4&5&6&7&8&9&10\\
\hline
$\dim M(1)^+_m$&1&0&1&1&3&3&6&7&12&14&22\\
\hline
\end{tabular}
\end{center}

\bigskip

In order to produce the second relation we need  the following lemma
whose proof is given in the Appendix.

\begin{lemma}\label{lemma:4.3} The vectors 
\begin{equation}\label{basis}
L(-1)M(1)^+_9,\quad L(-3)M(1)^+_7 \quad h(-1)^4_{-3}h(-1)^4\1,\quad
L(-2)^5\1
\end{equation}
span $M(1)^+_{10}.$
\end{lemma}

Now we can prove Proposition \ref{proposition:4.2}. 
First note that any weight $10$ vector is contained in $L(1,0)\oplus L(1,4)$
and is a linear combination of vectors of type $L(-n_1)\cdots L(-n_s)$,
$L(-m_1)\cdots L(-m_t)J$ where $n_1\geq \cdots \geq n_s\geq 2,$
$m_1\geq \cdots \geq m_t\geq 1,$ $\sum n_a=10$ and $\sum m_b+4=10.$  
{}From the proof of Theorem \ref{theorem:3.5} 
the images of this kind of vectors in
$A(M(1)^+)$ can be are expressed as linear combinations
of $\omega^i\,(i=0,1,\dots,5),\,\omega^i J\,(i=0,1,2,3)$.

By Proposition \ref{proposition:2.2} (i) and (iii),
$L(-1)M(1)^+_9,$ $L(-3)M(1)^+_7,$ $h(-1)^4_{-3}h(-1)^4\1$ are
congruent to vectors whose homogeneous components have weights less
than 10. Note that $\omega^i=L(-2)^i\1+$lower weight terms. 
Then it follows from 
Remark \ref{remark:3.6} and Proposition \ref{proposition:4.1}
that $L(-1)M(1)^+_9,$ $L(-3)M(1)^+_7,$ $h(-1)^4_{-3}h(-1)^4\1$ 
are congruent to vectors spanned by $\omega^i\,(i=0,1,\dots,4),\quad 
\omega^iJ\,(i=0,1,2).$ 
Thus by Lemma \ref{lemma:4.3} we see that
$M(1)^+_{10}+O(M(1)^+)$ is spanned by
$\omega^i\,(i=0,1,\dots,5),\omega^iJ\,(i=0,1,2).$
As a result we have
\[
\omega^3J = P(\omega)+Q(\omega)J
\]
where $\deg P \leq 5,\,\deg Q \leq 2$.  Evaluating this equation on
the top levels of the modules listed in Subsection \ref{subsection:3.4}
gives the desired result.

Now we can state our first main theorem.
\begin{theorem} \label{theorem:4.4}
We have the following algebra isomorphism:
\[
\C[x,y]/\langle P,Q\rangle\cong
A(M(1)^+)
\]
where
\[
P=(y+x-4x^2)(70y+908x^2-515x+27),
\quad
Q=(x-1)(x-\frac{1}{16})(x-\frac{9}{16})(y+x-4x^2).
\]
\end{theorem}

\begin{proof}
By Theorem \ref{theorem:3.5}, we have a surjective algebra homomorphism
\[
\begin{array}{cccc}
\varphi:&\C[x,y]&\longrightarrow&A(M(1)^+)\\
&x&\longmapsto&\omega\\
&y&\longmapsto&J.
\end{array}
\]
Let $K(x,y)\in \Ker\varphi$ and regard $K(x,y)$ as a polynomial in variable
$y.$ Note that $P(x,y)$ has degree 2 in $y.$ Using the division
algorithm we can write 
$K(x,y) = A(x,y)P(x,y)+R(x,y)$ where $A(x,y),R(x,y)\in \C[x,y]$ so that
$R(x,y)$ has degree 1 in $y.$ We can express $R(x,y)$ as 
$R(x,y) = B(x)(y+x-4x^2)+C(x)$. By Proposition \ref{proposition:4.1}
$P(x,y)\in \Ker\varphi.$ So we have
\begin{equation}\label{eqn:4.2}
B(\omega)(J+\omega-4\omega^2) +C(\omega) = 0.
\end{equation}
Evaluating the above equation on the top levels of 
modules $M(1,\lambda)$ yields
$C(\lambda^2/2) = 0$ since $J+\omega-4\omega^2=0$ on the top level of $M(1,\lambda)$
for all $\lambda\in\C^\times$. Thus $C(x)= 0$ as a polynomial.
Further evaluating equation (\ref{eqn:4.2}) on the top levels 
of $M(1)^-,\,M(1)(\theta)^\pm$ and noting that $J+\o-4\o^2\ne 0,$ we get
$B(1) = B(1/16)=B(9/16) = 0.$
This implies 
$(x-1)(x-1/16)(x-9/16)|B(x)$.
Thus we reach to
\[
K(x,y) = A(x,y)P(x,y)+D(x)Q(x,y)
\]
for some polynomial $D(x)$. Since $Q(x,y)$ lies in $\Ker\varphi$
already by Proposition \ref{proposition:4.2}, we conclude that
$ \Ker\varphi = \langle P(x,y),Q(x,y) \rangle.$ 
\end{proof}

\subsection{Classification of irreducible modules for $M(1)^+$}
\label{subsection:4.2}

Finally we can use $A(M(1)^+)$ whose algebra structure
was determined in the previous section to classify 
the irreducible modules for $M(1)^+.$

\begin{theorem} The set  
\[
\{M(1)^{\pm}, M(1)(\theta)^{\pm},
M(1,\l)\cong M(1,-\l), \l\in \C^\times\}
\]
gives a complete list of inequivalent irreducible $M(1)^+$
modules. Moreover, any irreducible admissible $M(1)^+$-module
is an ordinary module.
\end{theorem}
\begin{proof} Let $M=\oplus_{n\geq 0}M(n)$ be an irreducible
admissible $M(1)^+$-module with $M(0)\ne 0.$ Then $M(0)$ is an
irreducible $A(M(1)^+)$-module. Since $A(M(1)^+)$ is
commutative, $M(0)$ is one-dimensional. 
So both $\o$ and $J$ act as scalars $\alpha$ and $\beta$ on 
$M(0).$ From Theorem \ref{theorem:4.4} we have
$$
(\b+\a-4\a^2)(70\b+908\a^2-515\a+27)=0
$$
and
$$
(\a-1)(\a-\frac{1}{16})(\a-\frac{9}{16})(\b+\a-4\a^2)=0.
$$
If $\b+\a-4\a^2=0$ and $\a\ne 0$, then $M(0)$ is isomorphic to 
the top level of
$M(1,\sqrt{2\a})$ and $M$ is isomorphic to $M(1,\sqrt{2\a}).$
If $\b+\a-4\a^2=0$ and $\a=0$ then $M$ is isomorphic to $M(1)^+.$
Otherwise we have $(\a-1)(\a-1/16)(\a-9/16)=0$
and $70\b+908\a^2-515\a+27=0.$ One can easily verify that $M$ is
isomorphic to $M(1)^-,$ $M(1)(\theta)^+$ and $M(1)(\theta)^-$
when $\a=1,1/16$ and $9/16.$
\end{proof}

\section*{Appendix}
Here we give the details of a proof of Lemma \ref{lemma:4.3}.
First, we list bases of $M(1)^+_{7}$,$M(1)^+_{9}$ and
$M(1)^+_{10}$ which have dimension
$7,14$ and $22$ respectively.

\noindent
A basis of $M(1)^+_{7}$:
\[
\begin{array}{lll}
e_1= h(-6)h(-1)\1,&\quad e_2 = h(-5)h(-2)\1,&\quad e_3 = h(-4)h(-3)\1\\
e_4 = h(-4)h(-1)^3\1,&\quad e_5=h(-3)h(-2)h(-1)^2\1,&\quad e_6 = h(-2)^3h(-1)\1,\\
e_7 = h(-2)h(-1)^5\1.&&
\end{array}
\]

\noindent
A basis of $M(1)^+_{9}$:
\[
\begin{array}{ll}
f_1 = h(-8)h(-1)\1,& f_2= h(-7)h(-2)\1,\\
f_3= h(-6)h(-3)\1,&f_4 = h(-6)h(-1)^3\1,\\
f_5 = h(-5)h(-4)\1,& f_6=h(-5)h(-2)h(-1)^2\1,\\
f_7 = h(-4)h(-3)h(-1)^2\1,& f_8 = h(-4)h(-2)^2h(-1)\1,\\
f_9 = h(-4)h(-1)^5\1,&f_{10} =h(-3)^2h(-2)h(-1)\1,\\
f_{11}=h(-3)h(-2)^3\1, & f_{12}=h(-3)h(-2)h(-1)^4\1,\\
f_{13} = h(-2)^3h(-1)^3\1, & f_{14} = h(-2)h(-1)^7\1.
\end{array}
\]

\noindent
A basis of $M(1)^+_{10}$:
\[
\begin{array}{ll}
g_1 = h(-9)h(-1)\1,&\quad g_2 = h(-8)h(-2)\1,\\
g_3 = h(-7)h(-3)\1,&\quad g_4 = h(-7)h(-1)^3\1,\\
g_5 = h(-6)h(-4)\1,&\quad g_6 = h(-6)h(-2)h(-1)^2\1,\\
g_7 = h(-5)^2\1,&\quad g_8 = h(-5)h(-3)h(-1)^2\1,\\
g_9 = h(-5)h(-2)^2h(-1)\1,&\quad g_{10} = h(-5)h(-1)^5\1,\\
g_{11} = h(-4)^2h(-1)^2\1,&\quad g_{12} = h(-4)h(-3)h(-2)h(-1)\1,\\
g_{13}= h(-4)h(-2)^3\1,&\quad g_{14}=h(-4)h(-2)h(-1)^4\1,\\
g_{15}= h(-3)^3h(-1)\1,&\quad g_{16}=h(-3)^2h(-2)^2\1,\\
g_{17}=h(-3)^2h(-1)^4\1,&\quad g_{18}=h(-3)h(-2)^2h(-1)^3\1,\\
g_{19}=h(-3)h(-1)^7\1,&\quad g_{20}=h(-2)^4h(-1)^2\1,\\
g_{21}= h(-2)^2h(-1)^6\1,&\quad g_{22}=h(-1)^{10}\1.
\end{array}
\]

It is easy to see that
\begin{equation*}
\begin{split}
h&(-1)^4_{-3}h(-1)^4\1\\
=&96h(-9)h(-1)\1+ 144h(-7)h(-1)^3\1+144h(-6)h(-2)h(-1)^2\1\\
&+144h(-5)h(-3)h(-1)^2\1+72h(-4)^2h(-1)^2+48h(-5)h(-1)^5\1\\
&+96h(-4)h(-2)h(-1)^4\1+48h(-3)^2h(-1)^4\1+48h(-3)h(-2)^2h(-1)^3\1\\
&+4h(-3)h(-1)^7\1+6h(-2)^2h(-1)^6\1.
\end{split}
\end{equation*}
The tables below give the precise linear combinations of certain
vectors in terms of $g_i$ for $i=1,...,22.$ For example, 
$L(-1)f_1=8g_1+g_2.$ We know from the tables that the vectors in
(\ref{basis}) without $L(-2)^5\1$ 
span  a 21 dimensional subspace of $M(1)^+_{10}$
and none of these vectors involves the term $h(-1)^{10}\1$.
On the other hand, $L(-2)^5\1$ involves the term $h(-1)^{10}\1$.
Thus the vectors in (\ref{basis}) span $M(1)^+_{10},$ as expected.


\begin{center}
\begin{tabular}{|c||c|c|c|c|c|c|c|c|c|c|c|}
\hline
&$g_1$&$g_2$&$g_3$&$g_4$&$g_5$&$g_6$&$g_7$&$g_8$&$g_9$&$g_{10}$&$g
_{11}$\\
\hline
$L(-1)f_1$&8&1&0&0&0&0&0&0&0&0&0\\
\hline
$L(-1)f_2$&0&7&2&0&0&0&0&0&0&0&0\\
\hline
$L(-1)f_3$&0&0&6&0&3&0&0&0&0&0&0\\
\hline
$L(-1)f_4$&0&0&0&6&0&3&0&0&0&0&0\\
\hline
$L(-1)f_5$&0&0&0&0&5&0&4&0&0&0&0\\
\hline
$L(-1)f_6$&0&0&0&0&0&5&0&2&2&0&0\\
\hline
$L(-1)f_7$&0&0&0&0&0&0&0&4&0&0&3\\
\hline
$L(-1)f_8$&0&0&0&0&0&0&0&0&4&0&0\\
\hline
$L(-1)f_9$&0&0&0&0&0&0&0&0&0&4&0\\
\hline
$L(-1)f_{10}$&0&0&0&0&0&0&0&0&0&0&0\\
\hline
$L(-1)f_{11}$&0&0&0&0&0&0&0&0&0&0&0\\
\hline
$L(-1)f_{12}$&0&0&0&0&0&0&0&0&0&0&0\\
\hline
$L(-1)f_{13}$&0&0&0&0&0&0&0&0&0&0&0\\
\hline
$L(-1)f_{14}$&0&0&0&0&0&0&0&0&0&0&0\\
\hline
\hline
$L(-3)e_1$&6&0&0&0&1&1&0&0&0&0&0\\
\hline
$L(-3)e_2$&0&5&0&0&0&0&2&0&1&0&0\\
\hline
$L(-3)e_3$&0&0&4&0&3&0&0&0&0&0&0\\
\hline
$L(-3)e_4$&0&0&0&4&0&0&0&0&0&0&3\\
\hline
$L(-3)e_5$&0&0&0&0&0&3&0&2&0&0&0\\
\hline
$L(-3)e_6$&0&0&0&0&0&0&0&0&6&0&0\\
\hline
$L(-3)e_7$&0&0&0&0&0&0&0&0&0&2&0\\
\hline
\hline
$h(-1)^4_{-3}h(-1)^4\1$&96&0&0&144&0&144&0&144&0&48&72\\
\hline
\end{tabular}

\bigskip
Table A1
\end{center}

\begin{center}
\begin{tabular}{|c||c|c|c|c|c|c|c|c|c|c|c|c|c|c|c|c|c|c|c|c|c|c|}
\hline
&$g_{12}$&$g_{13}$
&$g_{14}$&$g_{15}$&$g_{16}$&$g_{17}$&$g_{18}$&$g_{19}$&$g_{20}$&$g_{
21}$&$g_{22}$\\
\hline
$L(-1)f_1$&0&0&0&0&0&0&0&0&0&0&0\\
\hline
$L(-1)f_2$&0&0&0&0&0&0&0&0&0&0&0\\
\hline
$L(-1)f_3$&0&0&0&0&0&0&0&0&0&0&0\\
\hline
$L(-1)f_4$&0&0&0&0&0&0&0&0&0&0&0\\
\hline
$L(-1)f_5$&0&0&0&0&0&0&0&0&0&0&0\\
\hline
$L(-1)f_6$&0&0&0&0&0&0&0&0&0&0&0\\
\hline
$L(-1)f_7$&2&0&0&0&0&0&0&0&0&0&0\\
\hline
$L(-1)f_8$&4&1&0&0&0&0&0&0&0&0&0\\
\hline
$L(-1)f_9$   &0&0&5&0&0&0&0&0&0&0&0\\
\hline
$L(-1)f_{10}$&6&0&0&2&1&0&0&0&0&0&0\\
\hline
$L(-1)f_{11}$&0&3&0&0&6&0&0&0&0&0&0\\
\hline
$L(-1)f_{12}$&0&0&3&0&0&2&4&0&0&0&0\\
\hline
$L(-1)f_{13}$&0&0&0&0&0&0&6&0&3&0&0\\
\hline
$L(-1)f_{14}$&0&0&0&0&0&0&0&2&0&7&0\\
\hline
\hline
$L(-3)e_1$&0&0&0&0&0&0&0&0&0&0&0\\
\hline
$L(-3)e_2$&0&0&0&0&0&0&0&0&0&0&0\\
\hline
$L(-3)e_3$&1&0&0&0&0&0&0&0&0&0&0\\
\hline
$L(-3)e_4$&0&0&1&0&0&0&0&0&0&0&0\\
\hline
$L(-3)e_5$&2&0&0&0&0&0&1&0&0&0&0\\
\hline
$L(-3)e_6$&0&1&0&0&0&0&0&0&1&0&0\\
\hline
$L(-3)e_7$&0&0&5&0&0&0&0&0&0&1&0\\
\hline
\hline
$h(-1)^4_{-3}h(-1)^4\1$&0&0&96&0&0&48&48&4&0&6&0\\
\hline
\end{tabular}

\bigskip
Table A2
\end{center}



\begin{thebibliography}{100000}
\bibitem[BFKNRV]{BFKNRV} R. Blumenhagen, M. Flohr, A. Kliem, W. Nahm,
A. Recknagel and R. Varnhagen, $W$-algebras with two and three 
generators, {\em Nucl. Phys.} {\bf B361} (1991), 255-289. 

\bibitem[B]{B86}
Borcherds,~R.: Vertex algebras, Kac-Moody
algebras, and the Monster, {\em Proc. Natl. Acad.
Sci. USA} 
\textbf{83}, 3068--3071(1986)

\bibitem[DVVV]{DVVV}
Dijkgraaf,~R., Vafa,~C., Verlinde,~E., Verlinde,~H.:
The operator algebra of orbifold models,
{\em Commu. Math. Phys.} \textbf{123}, 485--526(1989)

\bibitem[D]{D} C. Dong, Vertex algebras associated with 
even lattices, {\it J. Algebra} {\bf 160} (1993), 245-265.

\bibitem[DG]{DG}
Dong,~C., Griess,~R.L.,Jr.:
 Rank one lattice type vertex operator algebras and their
automorphism groups, {\em J. Algebra,} to appear, q-alg/9710017. 

\bibitem[DL]{DL} Dong,~C., Lepowsky,~J.:
 Generalized Vertex
Algebras and Relative Vertex Operators, Progress in Math. Vol. 112,
Birkh\"{a}user, Boston 1993.

\bibitem[DLM1]{DLM1} Dong,~C.,  Li,~H., Mason,~G.:
Regularity of rational vertex operator algebras, {\em Advances. in Math.}
{\bf 132} (1997), 148-166 

\bibitem[DLM2]{DLM2}
Dong,~C., Li,~H., Mason,~G.: 
Twisted representation of vertex operator algebras, 
{\em Math. Ann.} {\bf 310} (1998), 571-600.

\bibitem[DLi]{DLi}  Dong,~C., Lin,~Z.:
 Induced modules for vertex operator algebras, 
{\em Commu. Math. Phys.}
{\bf 179} (1996), 157-184.

\bibitem[DM1]{DM1} Dong,~C., Mason,~G.:
 Nonabelian orbifolds and boson-fermion 
correspondence,
{\em Commu. Math. Phys.} {\bf 163} (1994), 523-559.


\bibitem[DM2]{DM2} Dong,~C., Mason,~G.: 
 On quantum Galois theory, {\em Duke 
Math. J.} {\bf 86} (1997), 305-321. 

\bibitem[DN]{DN}  Dong,~C., Nagatomo,~K.:
 Representations of vertex operator algebra $V_L^+$ for rank one
lattice $L$, preprint.

\bibitem[FHL]{FHL}
Frenkel,~I.B., Huang,~Y., Lepowsky,~J.:
On axiomatic approach to vertex operator algebras and modules,
{\em Mem. Amer. Math. Soc.} \textbf{104} No.494(1993)


\bibitem[FLM]{FLM}
Frenkel,~I.B., Lepowsky,~J., Meurman,~A.:
Vertex operator algebras and the Monster,
Academic Press, 1988

\bibitem[FZ]{FZ}

Frenkel,~I.B, Zhu,~Y.:
Vertex operator algebras associated to
representations of affine and Virasoro
algebras, {\em Duke Math. J.} \textbf{66}, 123--168(1992)x

\bibitem[KL]{KL} M. Karel and H. Li, Certain generating subspaces for 
vertex operator algebras, preprint.

\bibitem[W]{W}
Wang,~W.:
Rationality of Virasoro vertex operator algebras,
{\em Duke Math. J.} \textbf{71}, IMRN No. 7, 197--211(1993)

\bibitem[Z]{Z}
Zhu,~Y.:
Modular invariance of characters of
vertex operator algebras,
{\em J.  AMS} \textbf{9}, 237--301(1996)
\end{thebibliography}
\end{document}